\documentclass[letterpaper,11pt]{article}
\usepackage{times}
\usepackage{empheq}
\usepackage{amsmath, amssymb}
\usepackage{amsthm}
\usepackage[T1]{fontenc}
\usepackage[english]{babel}
\usepackage[utf8]{inputenc}
\usepackage{hyperref}
\usepackage{mathrsfs}
\usepackage{mathabx}
\usepackage{color}
\usepackage{anysize}
\marginsize{2.5cm}{2.5cm}{2.5cm}{2.5cm}
\usepackage{url}
\usepackage{mathtools}
\usepackage[font={small,sc}]{caption}
\usepackage{ragged2e}
\usepackage[shortlabels]{enumitem}

\newtheorem{Theo}{Theorem}[section]

\AtBeginDocument{}

\newcommand{\Rd}{\mathbb{R}^{d}}

\newcommand{\R}{\mathbb{R}}

\begin{document}

\begin{Large}
\begin{center}
\textbf{Necessary and sufficient conditions for a family of continuous functions to form a Karhunen-Loève basis}
\end{center}
\end{Large}
\begin{center}
\large
\textsc{Ricardo CARRIZO VERGARA} \\
\normalsize
\textit{Université de Paris II Panthéon-Assas} \\
\textit{racarriz@uc.cl}
\end{center} 
\begin{center}
August 2022
\end{center}

The following fact is widely known and used in both practice and theory of stochastic processes: if $D \subset \Rd$ is a compact set, and if $K : D \times D  \to \R$ is a continuous positive-definite Kernel (also called a covariance Kernel), then $K$ has a decomposition in Fourier series of $L^{2}(D\times D)$ of the form
\begin{equation}
\label{Eq:KExpansion}
K = \sum_{n \in \mathbb{N}} \lambda_{n} f_{n} \otimes f_{n},
\end{equation}
where $(\lambda_{n})_{n \in \mathbb{N}} \in \ell^{1}(\mathbb{N})$ is a sequence of non-negative numbers, and $(f_{n})_{n \in \mathbb{N}}$ is an orthonormal basis of $L^{2}(D)$. The pair $(f_{n})_{n} , (\lambda_{n})_{n}$ consists in the eigenfunctions and eigenvalues of the associated covariance operator (Hilbert-Schmidt Theorem \cite[Theorem VI.16]{reed1980methods}), having
\begin{equation}
\label{Eq:fnLnEigen}
\lambda_{n}f_{n}(x) = \int_{D}K(x,y)f_{n}(y)dy.
\end{equation}
From equation \eqref{Eq:fnLnEigen} it is easy to verify that for $\lambda_{n} > 0$, the function $f_{n}$ is continuous. The cases in which $\lambda_{n} = 0$ do not really intervene in the decomposition \eqref{Eq:KExpansion}, hence the decomposition consists essentially of continuous functions. For simplicity, we will suppose $\lambda_{n} > 0 $ for every $n\in \mathbb{N}$ and that the family $(f_{n})_{n}$ is thus constituted of continuous functions (the family is hence not necessarily a \textit{complete} orthonormal system of $L^{2}(D)$). In addition, it is possible to prove that the convergence \eqref{Eq:KExpansion} is not only in the sense of $L^{2}(D\times D)$ but also uniform over $D\times D$ (Mercer's Theorem, which uses Dini's Theorem to prove the result \cite{mercer1909functions}). Those facts are the basis of the well-known and used Karhunen-Loève expansion of mean-square continuous stochastic processes over compact sets of the Euclidean space \cite{loeve1978probability,giambartolomei2015karhunen}.

The question we answer here is the one concerning the converse of such result: given an orthonormal system of $L^2(D)$ consisting of continuous functions $(f_{n})_{n}$, and given a sequence of strictly positive numbers $(\lambda_{n})_{n} \in \ell^{1}(\mathbb{N})$, does the limit
\begin{equation}
\label{Eq:LimitExpansion}
\lim_{n \to \infty} \sum_{j \leq n} \lambda_{j} f_{j} \otimes f_{j}
\end{equation}
exists uniformly over $D\times D$, defining thus a continuous covariance Kernel? Note that the limit \eqref{Eq:LimitExpansion} always exists in the $L^{2}(D\times D)$ sense, so the real question is if the resulting limit is continuous. By Dini's Theorem, using similar arguments as those used in the proof of Mercer's Theorem, one can conclude that if such limit is continuous, the convergence must be uniform.

The next result shows that the functions $(f_{n})_{n}$ and the coefficients $(\lambda_{n})_{n}$ must satisfy an extra criterion. Hence, not any pair basis-coefficients can be used to construct a Karhunen-Loève expansion.

\begin{Theo}
\label{Theo:KLbasisNecSuf}
Let $D \subset \Rd$ be a compact set. Let $(f_{n})_{n \in \mathbb{N}}$ be an orthonormal system of $L^{2}(D)$ consisting of continuous functions, and let $ ( \lambda_{n} )_{n \in \mathbb{N}} \in \ell^{1}(\mathbb{N})$ be a sequence of strictly positive numbers. Then, the limit \eqref{Eq:LimitExpansion} converges (uniformly) over $D\times D$ to a continuous covariance Kernel if and only if the sequence of functions
\begin{equation}
\label{Eq:DefVn}
v_{n} := \sum_{j \leq n} \lambda_{j} f_{j}^{2}, \quad n \in \mathbb{N},
\end{equation}
is equicontinuous over $D$. 
\end{Theo}

\textbf{Proof:} For the necessity, suppose $K$ is a continuous covariance kernel having the expansion \eqref{Eq:KExpansion}. Using Cauchy-Schwartz inequality, one concludes for every $x,y \in D$,
\begin{equation}
\label{Eq:BoundingIncrementsVn}
\begin{aligned}
|v_{n}(x)-v_{n}(y)| &= \left| \sum_{j \leq n} \lambda_{j}( f^{2}_{j}(x)-f^{2}_{j}(y) ) \right| \\
&\leq \sqrt{ \sum_{j \leq n} \lambda_{j}( f_{j}(x)-f_{j}(y) )^{2} }  \sqrt{ \sum_{j \leq n} \lambda_{j}( f_{j}(x)+f_{j}(y) )^{2} }  \\
&\leq \sqrt{ \sum_{j \in \mathbb{N}} \lambda_{j}( f_{j}(x)-f_{j}(y) )^{2} }  \sqrt{ 2\sum_{j \in \mathbb{N}} \lambda_{j}( f_{j}^{2}(x)+f_{j}^{2}(y) ) }  \\
&= \sqrt{ K(x,x) - 2K(x,y) + K(y,y)} \sqrt{2(K(x,x)+K(y,y))} \\
&\leq 2\sqrt{\| K \|_{\infty}}\sqrt{ K(x,x) - 2K(x,y) + K(y,y)}. 
\end{aligned}
\end{equation}
Since $K$ is uniformly continuous over the compact $D\times D$, for any given $\epsilon >0 $ there exists $\delta > 0 $ such that $|K(u,v)-K(x,y)| < \frac{\epsilon^{2}}{ 8\| K \|_{\infty}  }$ if $ |u-x|+|v-y| < \delta $. Hence, if $|x-y| < \delta $,
\begin{equation}
\label{Eq:KunifContinuous}
 |K(x,x) - 2K(x,y) + K(y,y)| \leq |K(x,x) - K(x,y)| + |K(y,y) - K(x,y)| < \frac{\epsilon^{2}}{ 4\| K \|_{\infty}  }.
\end{equation}
Using this in inequality \eqref{Eq:BoundingIncrementsVn}, we obtain that if $|x-y| < \delta$ then
\begin{equation}
\label{Eq:VnEquicontinuous}
|v_{n}(x)-v_{n}(y)| \leq \epsilon, \quad \forall n \in \mathbb{N},
\end{equation}
hence the family $(v_{n})_{n}$ is indeed equicontinuous.

For the sufficiency, suppose the sequence $(v_{n})_{n}$ is equicontinuous over $D$. For $n,m \in \mathbb{N}$, using the orthonormality of the family $(f_{n})_{n}$ and the fact that $\sum_{n \in \mathbb{N}} \lambda_{n} < \infty$, we have
\begin{equation}
\int_{D}| v_{n} - v_{m} |dx = \int_{D} \sum_{j=n \wedge m + 1}^{n \vee m} \lambda_{j} f_{j}^{2} dx = \sum_{j=n \wedge m + 1}^{n \vee m} \lambda_{j} \xrightarrow[n,m \to \infty]{} 0.
\end{equation}
Hence $(v_{n})_{n}$ is a Cauchy sequence of $L^{1}(D)$, hence it converges by completeness to a function $v \in L^{1}(D)$. Now, the converges in $L^{1}(D)$ implies the existence of a sub-sequence $(v_{n_{k}})_{k}$ which converges almost-everywhere in $D$ to a function which must coincide with $v$ almost-everywhere. Let $\tilde{D} \subset D$ be the full-measure set where the sub-sequence $(v_{n_{k}})_{k}$ converges point-wise to $v$. Hence, for every $x \in \tilde{D}$, the sequence $(v_{n}(x))_{n}$ has a sub-sequence which converges to $v(x)$. But since the sequence $v_{n}(x)$ consists in the partial sums of a series of positive numbers, then the convergence of a sub-sequence implies the convergence of the entire sequence $(v_{n}(x))_{n}$. We conclude hence that $v_{n}$ converges to $v$ point-wise over $\tilde{D}$. Since $\tilde{D}$ has full-measure, it must be dense in $D$. Since the family $(v_{n})_{n}$ is equicontinuous and it converges point-wise to a function $v$ in a dense subset of $D$, then the convergence is everywhere in $D$ and the function $v$ must be continuous \cite[Theorem I.26]{reed1980methods}. In addition, since the sequence $(v_{n})_{n}$ is monotonically increasing, by Dini's theorem the convergence must be uniform.

The uniform convergence of $(v_{n})_{n}$ and the Cauchy-Schwarz inequality allow us to conclude
\small
\begin{equation}
\label{Eq:UniformCauchyConvExpansionK}
\begin{aligned}
\sup_{(x,y) \in D\times D} \left| \sum_{j= n\wedge m   + 1 }^{n \vee m } \lambda_{j}f_{j}(x)f_{j}(y) \right| &\leq \sup_{(x,y) \in D\times D} \sqrt{  \sum_{j= n\wedge m   + 1 }^{n \vee m } \lambda_{j}f_{j}^{2}(x)} \sqrt{  \sum_{j= n\wedge m   + 1 }^{n \vee m } \lambda_{j}f_{j}^{2}(y)} \\
&= \| v_{n} - v_{m} \|_{\infty} \xrightarrow[n,m \to \infty]{} 0
\end{aligned}.
\end{equation}
\normalsize
From the Cauchy-sequence convergence criterion, the series $\sum_{j \in \mathbb{N}} \lambda_{j} f_{j}\otimes f_{j}$ converges uniformly over $D\times D$ to some function $K$ which must hence be continuous, and by construction, it is also a positive-definite Kernel. $\blacksquare$

\bibliography{mibib}
\bibliographystyle{abbrv}

\end{document}